\documentclass[a4paper]{amsart}
\usepackage[T1]{fontenc}
\usepackage{amscd}
\usepackage{pstricks}
\usepackage{hyperref}

\title[Real points of moduli schemes of vector bundles]{Real points of coarse moduli schemes of vector bundles on a real algebraic curve}

\author{Florent Schaffhauser}
\address{Department of Mathematics, University of Los Andes, Bogota, Colombia.}
\email{florent@uniandes.edu.co}

%{Algebraic moduli problems, moduli of vector bundles; Special connections and metrics on vector bundles (Hermite-Einstein-Yang-Mills)}

\newtheorem{theorem}{Theorem}[section]

\newtheorem{proposition}[theorem]{Proposition}
\newtheorem{corollary}[theorem]{Corollary}
\newtheorem{definition}[theorem]{Definition}

\newtheorem*{thm*}{Theorem}

\theoremstyle{definition}
\newtheorem*{ack}{Acknowledgments}

\newcommand{\C}{\mathbb{C}}
\newcommand{\R}{\mathbb{R}}
\renewcommand{\H}{\mathbb{H}}
\newcommand{\Z}{\mathbb{Z}}

\newcommand{\Id}{\mathrm{Id}}

\newcommand{\gr}{\mathrm{gr}}
\newcommand{\tr}{\mathrm{tr} }
\newcommand{\rk}{\mathrm{rk}}
\newcommand{\Gal}{\mathrm{Gal}}

\newcommand{\Spec}{\mathrm{Spec}\, }
\newcommand{\Pic}{\mathrm{Pic}}

\newcommand{\Stab}{\mathrm{Stab}}

\newcommand{\End}{\mathrm{End}}

\newcommand{\calF}{\mathcal{F}}

\newcommand{\calE}{\mathcal{E}}
\newcommand{\cE}{\mathcal{E}}
\newcommand{\cF}{\mathcal{F}}

\newcommand{\calO}{\mathcal{O}}
\newcommand{\calL}{\mathcal{L}}
\newcommand{\calI}{\mathcal{I}}
\newcommand{\cL}{\mathcal{L}}

\newcommand{\ov}[1]{\overline{#1}}
\newcommand{\os}[1]{\overline{\sigma^*#1}}

\renewcommand{\phi}{\varphi}
\renewcommand{\mod}[1]{\mathrm{mod\ #1}}
\renewcommand{\deg}{\mathrm{deg}}

\renewcommand{\mod}{\mathrm{mod\ }}
\newcommand{\w}{\omega}
\newcommand{\sig}{\sigma}
\newcommand{\sigt}{\widetilde{\sigma}}
\newcommand{\Ms}{M^{\sig}}

\newcommand{\conn}{\mathcal{A}_{E}}

\newcommand{\Mod}{\mathcal{M}^{\, r,d}_{M,\sig}}
\newcommand{\Mods}{\mathcal{N}^{\, r,d}_{M,\sig}}
\newcommand{\sectionsofE}{\Omega^0(M;E)}
\newcommand{\zerooneforms}{\Omega^{0,1}(M;E)}
\newcommand{\onezeroforms}{\Omega^{1,0}(M;E)}
\newcommand{\oneforms}{\Omega^{1}(M;E)}
\newcommand{\kforms}{\Omega^{k}(M;E)}

\newcommand{\asigt}{\alpha_{\sigt}}

\newcommand{\gaugegp}{\mathcal{G}_{E}}
\newcommand{\cxgaugegp}{\mathcal{G}^{\, \C}_{E}}
\newcommand{\realgaugegp}{\mathcal{G}_{E}^{\, \sigt}}
\newcommand{\antiHermoneforms}{\Omega^1(M;\mathfrak{u}(E))}
\newcommand{\antiHermtwoforms}{\Omega^2(M;\mathfrak{u}(E))}

\newcommand{\fibre}{F^{-1}\big(\{*i2\pi\frac{d}{r}\, \mathrm{Id}_E\}\big)}
\newcommand{\asigtp}{\alpha_{\sigt'}}
\newcommand{\phisigt}{\phi_{\sigt}}
\newcommand{\phisigtp}{\phi_{\sigt'}}
\newcommand{\tausigt}{\tau_{\sigt}}
\newcommand{\tausigtp}{\tau_{\sigt'}}
\newcommand{\Bun}{\mathcal{B}un_{ss,\mu}}
\newcommand{\BunR}{\mathcal{B}un_{ss,\mu}^{\ \R}}
\newcommand{\BunH}{\mathcal{B}un_{ss,\mu}^{\ \H}}
\newcommand{\Ker}{\mathrm{Ker}\,}
\renewcommand{\Im}{\mathrm{Im}\,}

\newtheorem{example}[theorem]{Example}

\begin{document}

\begin{abstract}
We examine a moduli problem for real and quaternionic vector bundles on a smooth complex projective curve with a fixed real structure, and we give a gauge-theoretic construction of moduli spaces for semi-stable such bundles with fixed topological type. These spaces embed onto connected subsets of real points inside a complex projective variety. We relate our point of view to previous work by Biswas, Huisman and Hurtubise (\cite{BHH}), and we use this to study the $\Gal(\C/\R)$-action $[\calE] \mapsto [\os{\calE}]$ on moduli varieties of \textit{stable} holomorphic bundles on a complex curve with given real structure $\sigma$. We show in particular a Harnack-type theorem, bounding the number of connected components of the fixed-point set of that action by $2^g +1$, where $g$ is the genus of the curve. In fact, taking into account all the topological invariants of $\sigma$, we give an exact count of the number of connected components, thus generalising to rank $r > 1$ the results of Gross and Harris on the Picard scheme of a real algebraic curve (\cite{GH}).
\end{abstract}

\maketitle

\tableofcontents

\section{Introduction}

Let $(X,  \calO_X)$ be a scheme of finite type over $\R$ such that $M := X(\C)$ is an irreducible, smooth, complex projective curve. We shall call such a scheme a real algebraic curve. The complex conjugation of $\C$ induces a continuous action of the group $G := \Gal(\C/\R)$ on the complexified scheme $$(X_{\C} := X \times_{\Spec \R} \Spec \C,\,  \calO_{X_{\C}} :=  \calO_X \otimes_{\R} \C)\, .$$ $M$ is the set of closed points of $X_{\C}$, so the Galois action on $X_{\C}$ induces an involution $\sig$  of $M$, whose tangent map is $\C$-antilinear. In differential geometric terms, $M$ is a compact, connected Riemann surface (with an integral K\"ahler metric, which we take to be of unit volume), and $\sigma$ is an anti-holomorphic, involutive isometry of $M$. The quotient space $M/\sig$ is the set of closed points of $X = X_{\C}/G$.  We denote $p$ the projection $p: X_{\C} \longrightarrow X=X_{\C}/G$. Then $G$ acts on the sheaf $p_*O_{X_{\C}}$, and one has $O_X \simeq (p_*O_{X_{\C}})^G$. Similarly, if $O_M$ denotes the sheaf of holomorphic functions on $M$, $G$ acts on $p_*O_M$ by $$(\sigma\cdot f)(x) = \ov{f\big(\sig(x)\big)},$$ where $f$ is any holomorphic function defined on a $\sigma$-invariant open subset of $M$ (in particular, $f$ is real-valued on $\Ms$, the fixed-point set of $\sigma$ in $M$). In differential geometric terms, the ringed space $(M/\sig,  (p_*\calO_M)^{G})$ is called a Klein surface (\cite{AG}). The boundary of $M/\sig$ is diffeomorphic to $X(\R)$, the set of real points of $X$. In particular, it might be empty. Topologically, $M^{\sig} = X(\R)$ is a disjoint union of $k$ circles embedded in $M=X(\C)$. By Harnack's theorem, one has $0 \leq k \leq g+1$, where $g$ is the genus of $M$. We note that $M/\sig$ topologically is a compact connected real surface which, necessarily, is either non-orientable or has non-empty boundary (it can be both, but orientable surfaces without boundary are excluded). In the figure below, we denote $|X|$ the set of closed points of a scheme $X$, and $\mathrm{Proj}\ R$ the homogeneous spectrum of a graded ring $R$.

\begin{figure}[ht]
\centerline{
\begin{pspicture}(0,-1)(7,2)
\pscircle(0,0){2}
\psarc(0,4){4.47213}{-116.56505}{-63.43495}
\rput(-2.5,0){$\sig \updownarrow$}
\rput(3,0){$M^{\sig}=X(\R)$}
\psarc[linestyle=dashed,linewidth=0.005](0,-4){4.47213}{63.43495}{116.56505}
\rput(0,-2.5){$M=X(\C)=|X_{\C}|$}
\psarc(6.5,0){2}{0}{180}
\psarc(6.5,4){4.47213}{-116.56505}{-63.43495}
\psarc[linestyle=dashed,linewidth=.005](6.5,-4){4.47213}{63.43495}{116.56505}
\rput(9,0){$X(\R)$}
\rput(7,-2.5){$M/\sig = |X|=|X_{\C}|/G$}
\end{pspicture}
}
\vspace{1.5cm}
\caption{$X=\mathbb{P}^1_{\R} =  \mathrm{Proj}\ (\R[X_0,X_1])$}
\end{figure}

\noindent For any $r \geq 1$ and any integer $d$, we denote $\Mod$ the coarse moduli scheme parametrising $S$-equivalence classes of semi-stable vector bundles of rank $r$ and degree $d$ on $M=X(\C)$, and $\Mods$ the open sub-scheme of $\Mod$ parametrising isomorphism classes of stable bundles of rank $r$ and degree $d$ on $M$. We recall that $\Mods=\Mod$ if, and only if, $r\wedge d=1$. As $X$ is defined over $\R$, so are the moduli schemes $\Mod$ and $\Mods$, and, in the present paper, we are interested in the topology of the set $\Mods(\R)$ as a subset of the complex variety $\Mods(\C)$, endowed with its complex topology. We recall that, for $g \geq 2$, $\Mods(\C)$ is a smooth, connected, complex quasi-projective variety of dimension $r^2(g-1) + 1$ (in particular, the dimension does not depend on $d$). Our main result is a precise count of the connected components of $\Mods(\R)$, which nicely generalises the results of Gross and Harris on the Picard scheme of a real algebraic curve (\cite{GH}). Note that in the $r=1$ case, more is known~: any given connected component of $\Pic_X^d(\R) = \mathcal{M}_{M,\sig}^{1,d}(\R)$ is a real $g$-dimensional torus $\R^g / \Z^g$, $g$ being the genus of $X$.

To be able to state our result, let us recall the topological classification of real algebraic curves, first obtained by Felix Klein (\cite{Klein1}). Given an algebraic curve $X$ defined over $\R$, let us denote $g(X)$ the genus of $X$, $k(X)$ the number of connected components of $X(\R)$, and $a(X)$ the number defined by
\begin{eqnarray*}
& & a(X) = 0\ \mathrm{if}\ X(\C) \setminus X(\R)\ \mathrm{is\ not\ connected,}\\
& & a(X) = 1\ \mathrm{if}\ X(\C) \setminus X(\R)\ \mathrm{is\ connected.}
\end{eqnarray*}

\noindent Equivalently, $a(X) = 0$ if $X(\C) / G$ is orientable, and $a(X) = 1$ if $X(\C) / G$ is non-orientable; for that reason is $a(X)$ sometimes called the orientability index of $X$. Given two real algebraic curves $X$ and $X'$, Klein's classification theorem says that there exists a Galois-equivariant homeomorphism between $X(\C)$ and $X'(\C)$ if, and only if, $$g(X) = g(X'),\ k(X) = k(X'),\ \mathrm{and}\ a(X) = a(X').$$ The classification of real compact connected surfaces shows that this is in fact equivalent to $X(\C) / G$ and $X'(\C) / G$ being homeomorphic. Moreover, one has~:

\begin{itemize}

\item $0 \leq k(X) \leq g(X) + 1$ (Harnack's theorem),

\item if $k(X) = 0$ then $a(X) = 1$, and if $k(X) = g+1$ then $a(X) = 0$,

\item if $a(X) = 0$, then $k(X) \equiv (g+1)\ (\mod{2})$,

\end{itemize}

\noindent and Klein proved that all triples $(g(X), k(X), a(X))$ satisfying the conditions above occur for some real algebraic curve $X$. Our main result then is as follows.

\begin{theorem}\label{main_result}
Let $X \leftrightarrow (M,\sigma)$ be a real algebraic curve of topological type $(g,k,a)$, $g\geq 2$, and let $\Mods$ be the coarse moduli scheme parametrising isomorphism classes of stable holomorphic vector bundles of rank $r$ and degree $d$ on $X(\C)$. We consider the real structure induced on the smooth, connected, complex quasi-projective variety $\Mods(\C)$ by the functor $\cE \longmapsto \os{\cE}$.

\begin{enumerate}

\item Assume that $k > 0$.

\begin{enumerate}

\item If $r \equiv 1\ (\mod{2})$, then $\Mods(\R)$ is non-empty and has $2^{k-1}$ connected components. For fixed $r$ and $d$, these connected components are pairwise homeomorphic.

\item If $r \equiv 0\ (\mod{2})$ and $d \equiv 1 (\mod{2})$, then $\Mods(\R)$ is non-empty and has $2^{k-1}$ connected components.

\item If $r \equiv 0\ (\mod{2})$ and $d \equiv 0\ (\mod{2})$, then $\Mods(\R)$ is non-empty and has $2^{k-1} + 1$ connected components. 

\end{enumerate}

\item Assume that $k=0$.

\begin{enumerate}

\item If $r(g-1) \equiv 0\ (\mod{2})$ and $d \equiv 0\ (\mod{2})$, then $\Mods(\R)$ is non-empty and has $2$ connected components. For fixed $r$ and $d$, and if $g \equiv 1 (\mod{2})$, these two connected components are homeomorphic.

\item If $r(g-1) \equiv 0\ (\mod{2})$ and $d \equiv 1\ (\mod{2})$, then $\Mods(\R)$ is empty.

\item If $r(g-1) \equiv 1\ (\mod{2})$, then $\Mods(\R)$ is non-empty and has $1$ connected component.

\end{enumerate}

\end{enumerate}

\end{theorem}

\noindent As a corollary to the above result combined with Harnack's theorem, we obtain an upper bound, depending only on the genus of $X$, on the number of connected components of $\Mods(\R)$.

\begin{corollary}
Given a real algebraic curve $X$ of genus $g$, the number of connected components of $\Mods(\R)$ is lower than $2^{g} +1$.
\end{corollary}

\noindent We recall that, when $g$, $r$ and $d$ remain fixed while $k$ or $a$ (the real structure of $X$) changes, $\Mods(\C)$ stays the same topologically, but the topology of the connected components of $\Mods(\R)$ may change (see Section 6.3 of \cite{BHH}).
\begin{ack}
It is a pleasure to thank Luis \'Alvarez-C\'onsul, Thomas Baird, Mario Garc\'ia-Fern\'an\-dez, \'Oscar Garc\'ia-Prada, Olivier Guichard, Johannes Huisman, Arturo Prat-Waldron, Jan Swoboda, and Richard Wentworth, for discussions upon the topics dealt with in the present paper. Special thanks go to Melissa Liu for her very careful reading of the paper and her precious comments. Thanks also to the \textit{Max Planck Society}, for supporting my stay at the \textit{Max Planck Institute for Mathematics} in Bonn, where the first version of this paper was written. Finally, I would like to thank the referee for pointing out a mistake in the first version of the paper, and for many valuable comments on the notion of stability studied in the present work.
\end{ack}

\section{Vector bundles on a real algebraic curve}\label{vb_over_real_curves}

\subsection{Real and quaternionic vector bundles}

Our motivation to study real points of coarse moduli schemes of vector bundles comes from an attempt at formulating a moduli problem for vector bundles on a real algebraic curve. The fundamental tool in that endeavour is the notion of \textit{real space}, which is due to Atiyah (\cite{Atiyah_real_bundles}). In the special context of real algebraic geometry (from which the notion of real space originates), Atiyah's observation is that the category of real algebraic vector bundles on a real algebraic curve $X$ is equivalent to the category of holomorphic vector bundles on $M=X(\C)$, endowed with a $\C$-antilinear involution covering the natural involution $\sig$ of $M$.

\begin{definition}[Real vector bundles]
A \textbf{real vector bundle} on $(M,\sig)$ is a pair $(\calE,\sigt)$ where $\calE$ is a holomorphic vector bundle and $\sigt:\calE \longrightarrow \calE$ is a map satisfying the following conditions~:

\begin{enumerate}

\item the diagram

\begin{equation*}
\begin{CD}
\calE 	@>\sigt>> 	\calE \\
@VVV 				@VVV \\
M 		@>\sig>> 		M
\end{CD}
\end{equation*}

\noindent is commutative,

\item the map $\sigt$ is $\C$-antilinear,

\item $\sigt^2 = \Id_{\calE}$.

\end{enumerate}

\end{definition}

\noindent We shall refer to the map $\sigt$ as the \textbf{real structure} of $\calE$. It induces a $\C$-linear isomorphism $\phi$, covering the identity of $M$, between $\os{\calE}$ and $\calE$. This isomorphism satisfies $\os{\phi} = \phi^{-1}$ and, as a matter of fact, giving a $\C$-linear, invertible morphism $\phi: \os{\calE} \to \calE$ satisfying $\os{\phi} = \phi^{-1}$, is equivalent to giving a $\C$-antilinear, invertible map $\sigt:\calE\longrightarrow\calE$ covering $\sigma$ and squaring to the identity. A homomorphism of real holomorphic vector bundles is a homomorphism of holomorphic vector bundles (covering the identity map of $M$ and) commuting to the respective real structures. One may observe that the two equivalent categories of real algebraic vector bundles on $X$ and real holomorphic vector bundles (in the sense of Atiyah) on $M=X(\C)$, are equivalent to a third one, namely the category of dianalytic vector bundles (complex vector bundles which admit an atlas whose transition maps are either holomorphic or anti-holomorphic) on the Klein (=dianalytic) surface $M / \sig$, as shown in \cite{Sch_GD}. In that sense, and provided that one accepts to work in the dianalytic category, studying vector bundles on a compact surface which is either non-orientable or has non-empty boundary is equivalent to studying vector bundles on a real algebraic curve. At any rate, what is starting to shape here is that real holomorphic vector bundles define real points of the moduli schemes $\Mods$. But, as it turns out, there might be real points of a slightly different type, due to the presence of non-trivial automorphisms for stable vector bundles on $M$. Indeed, a real vector bundle $\calE$ on $M$ certainly is self-conjugate (meaning that $\os{\calE} \simeq \calE$), but the converse is not true, even if $\calE$ only has scalar automorphisms, and this leads to the notion of quaternionic vectors bundles (or symplectic vectors bundles, as they are called in \cite{Dupont} and in \cite{Hartshorne}).

\begin{definition}[Quaternionic vector bundles]
A \textbf{quaternionic vector bundle} on $(M,\sig)$ is a pair $(\calE,\sigt)$ where $\calE$ is a holomorphic vector bundle and $\sigt:\calE \longrightarrow \calE$ is a map satisfying the following conditions~:

\begin{enumerate}

\item the diagram

\begin{equation*}
\begin{CD}
\calE 	@>\sigt>> 	\calE \\
@VVV 				@VVV \\
M 		@>\sig>> 		M
\end{CD}
\end{equation*}

\noindent is commutative,

\item the map $\sigt$ is $\C$-antilinear,

\item $\sigt^2 = - \Id_{\calE}$.

\end{enumerate}

\end{definition}

\noindent We shall refer to the map $\sigt$ as the \textbf{quaternionic structure} of $\calE$. It induces a $\C$-linear isomorphism $\phi$, covering the identity of $M$, between $\os{\calE}$ and $\calE$. This isomorphism satisfies $\os{\phi} = - \phi^{-1}$ and, as a matter of fact, giving a $\C$-linear, invertible morphism $\phi: \os{\calE} \to \calE$ satisfying $\os{\phi} = - \phi^{-1}$, is equivalent to giving a $\C$-antilinear, invertible map $\sigt:\calE\longrightarrow\calE$ covering $\sigma$ and squaring to minus the identity. A homomorphism of quaternionic holomorphic vector bundles is a homomorphism of holomorphic vector bundles (covering the identity map of $M$ and) commuting to the respective quaternionic structures. One may observe here that, when $\Ms = \emptyset$, the complex rank of a quaternionic bundle is allowed to be odd, while when $\Ms\not=\emptyset$, it must be even, for the fibres of $E|_{\Ms} \to \Ms$ are left modules over the field of quaternions. In the present work, stability always means slope stability (see Definition \ref{stability_def}).

\begin{proposition}[\cite{BHH}]\label{stable_and_Galois_invariant}
Assume that $\calE$ is a stable holomorphic bundle on $M$, and that $\os{\calE} \simeq \calE$. Then $\calE$ is either real or quaternionic, and it cannot be both.
\end{proposition}

\begin{proof}
We recall that a stable bundle only has a scalar automorphisms, because its endomorphism ring is a field (a non-zero morphism between stable bundles of equal slope is an isomorphism) that contains $\C$ as a sub-field (the sub-field of scalar endomorphisms), and its elements are algebraic over $\C$ by the Cayley-Hamilton theorem, so they are contained in $\C$. We then proceed with the proof of the Proposition. A $\C$-linear isomorphism $\phi: \os{\calE} \overset{\simeq}{\longrightarrow} \calE$ covering $\Id_M$ is the same as a $\C$-antilinear map $\sigt: \calE \to \calE$ covering $\sigma$. As $\sig^2 = Id_M$, the map $\sigt^2$ is a $\C$-linear map covering $\Id_M$. Since $\calE$ only has scalar automorphisms, this implies that $\sigt^2 = \lambda\in\C^*$. Replacing $\sigt$ with $\sigt/\sqrt{|\lambda|}$ if necessary, we may assume that $|\lambda| = 1$. Moreover, $\lambda \sigt = (\sigt^2) \sigt = \sigt (\sigt^2) = \sigt(\lambda \cdot) = \ov{\lambda} \sigt$, so $\lambda = \ov{\lambda}$. As a consequence, $\lambda = \pm 1$, making $\calE$ real or quaternionic. If $\sigt'$ is another $\C$-antilinear map covering $\sigma$, then, as $\calE$ only has scalar automorphisms, $\sigt'\circ \sigt = \nu \in\C^*$, so $$(\sigt')^2 (\sigt)^2 = \sigt' (\nu\cdot) \sigt = \ov{\nu}\, \sigt'\sigt = |\nu|^2 \Id_{\calE},$$ with $|\nu|^2>0$. Therefore, $\sigt$ and $\sigt'$ are either both real or both quaternionic.
\end{proof}

\noindent As the Galois action on $\Mods(\C)$ is induced by the functor $\calE \mapsto \os{\calE}$, which preserves rank, degree, and slope stability of a holomorphic vector bundle, we see that real points of $\Mods$ may consist of real and quaternionic vector bundles alike. The precise situation will become clearer after we have identified the connected components of $\Mods(\R)$.

\subsection{Topological classification}

Any attempt at finding moduli for real and qua\-ternionic bundles on $X(\C)$ begins with the determination of some discrete invariants specifying a topological (or smooth) type for those bundles. Such a classification result was obtained by Biswas, Hurtubise and Huisman in \cite{BHH} (Propositions 4.1, 4.2, and 4.3). We formulate their result in the additional presence of a smooth Hermitian metric on the vector bundles that we consider. In this context, a real or quaternionic structure $\sigt$ on a Hermitian, smooth complex vector bundle $E$ is assumed to be an isometry.

\begin{theorem}[\cite{BHH}]\label{top_type}
One has~:
\begin{itemize}

\item For real bundles~:

\begin{itemize}

\item if $\Ms = \emptyset$, then real Hermitian bundles on $(M,\sig)$ are topologically classified by their rank and degree. 
It is necessary and sufficient for a real Hermitian bundle of rank $r$ and degree $d$ to exist that 
$$
d \equiv 0\ (\mod{2}).
$$

\item if $\Ms \not= \emptyset$ and $(E,\sigt)$ is real, then $(E^{\sigt} \to M^{\sig})$ is a real vector bundle in the ordinary sense, on the disjoint union 
$$
\Ms = \gamma_1 \sqcup \cdots \sqcup \gamma_n
$$ 
of at most $(g+1)$ circles, and we denote $$w^{(j)} := w_1(E^{\sigt}\big|_{\gamma_j}) \in H^1(S^1; \Z / 2\Z) \simeq \Z / 2\Z$$ the first Stiefel-Whitney class of 
$E^{\sigt}\to \Ms$  restricted to $\gamma_j$.\\ Then real Hermitian bundles on $(M,\sig)$ are topologically classified by their 
rank, their degree, and the sequence  $(w^{(1)}, \cdots, w^{(n)})$. 
It is necessary and sufficient for a real Hermitian bundle with given invariants $r$, $d$ and 
$(w^{(1)}, \cdots, w^{(n)})$ to exist that 
$$
w^{(1)} + \cdots + w^{(n)} \equiv d\ (\mod 2).
$$
\end{itemize}

\item For quaternionic bundles~:\\ Quaternionic Hermitian bundles on $(M,\sig)$ are topologically classified by their rank and degree. 
It is necessary and sufficient for a topological quaternionic bundle of rank $r$ and degree $d$ to exist that 
$$
d + r(g-1) \equiv 0\ (\mod 2).
$$

\end{itemize}

\end{theorem} 

\subsection{Stability}\label{stability}

As a next step into the moduli problem for algebraic vector bundles on a real algebraic curve, it certainly is necessary to have a notion of stability at our disposal in order to proceed. In the context of vector bundles on a curve, slope stability probably is the obvious choice, but it is perhaps not so clear whether one should test that condition over all sub-bundles of $\calE$, or over real sub-bundles only (that is, over sub-bundles on which $\sigt$ induces a real structure). This is an important matter because, as it turns out, different choices at this stage lead to different answers later. We therefore devote some time to analysing the different notions of stability for real and quaternionic bundles, and comparing them.

The slope of a non-zero holomorphic vector bundle $\cE$ is the quotient $$\mu(\cE) := \frac{\deg\,\cE}{\rk\,\cE}$$ of its degree by its rank.

\begin{definition}[Stability conditions for real and quaternionic bundles]\label{stability_def}
Let $(\cE,\sigt)$ be a real (resp. quaternionic) holomorphic vector bundle on $(M,\sigma)$. We call a sub-bundle of $\cE$ non-trivial if it is distinct from $\{0\}$ and from $\cE$. Then $(\cE,\sigt)$ is said to be 
\begin{enumerate}
\item \textbf{stable} if, for any non-trivial $\sigt$-invariant sub-bundle $\cF \subset \cE$, the slope stability condition $$\mu(\cF) < \mu(\cE)$$ is satisfied.
\item \textbf{semi-stable} if, for any non-trivial  $\sigt$-invariant sub-bundle $\cF \subset \cE$, one has $$\mu(\cF) \leq \mu(\cE).$$
\item \textbf{geometrically stable} if the underlying holomorphic bundle $\cE$ is stable, that is, if, for any non-trivial sub-bundle $\cF \subset \cE$,  one has $$\mu(\cF) < \mu(\cE).$$
\item \textbf{geometrically semi-stable}, if the underlying holomorphic bundle $\cE$ is semi-stable, that is, if for any non-trivial sub-bundle $\cF \subset \cE$,  one has $$\mu(\cF) \leq \mu(\cE).$$
\end{enumerate}
\end{definition}

\noindent Various comments are in order. First, we notice that, when $(\cE,\sigt)$ is real, it is of the form $\calE=\mathbb{E}(\C)$ for some algebraic vector bundle $\mathbb{E}\longrightarrow X$ defined over the reals, and geometric stability means stability of the bundle $\mathbb{E}(\C)$ of geometric points of $\mathbb{E}$. Geometric stability is, for instance, the notion of stability chosen in \cite{HN} (Section 1.1, page 217). Second, we see that (3) $\Rightarrow$ (1), and (4) $\Rightarrow$ (2). We prove below that (2) $\Rightarrow$ (4), but (1) $\not\Rightarrow$ (3).

\begin{proposition}\label{semi-stability_and_geometric_semi-stability}
Let $(\cE,\sigt)$ be a semi-stable real (resp. quaternionic) vector bundle on $(M,\sig)$. Then $(\cE,\sigt)$ is geometrically semi-stable.
\end{proposition}

\begin{proof}
Let $\phi: \os{\cE} \longrightarrow \cE$ be the isomorphism determined by the real (resp. quaternionic) structure on $\cE$. Assume that $(\cE,\sigt)$ is not geometrically semi-stable, and let $\calF$ be the destabilising bundle of $\calE$ (the unique maximal rank bundle among sub-bundles of $\calE$ the slope of which is maximal). Then $\phi(\os{\calF})$ and $\calF$ are sub-bundles of $\calE$ which have the same rank and degree. By unicity of $\calF$, one has $\phi(\os{\cF}) = \cF$. So $\calF$ is $\sigt$-invariant, and therefore $\mu(\calF) \leq \mu(\calE)$, which contradicts the assumption that $\calF$ is the destabilising bundle for $\calE$. 
\end{proof}

\noindent Proposition \ref{semi-stability_and_geometric_semi-stability} is actually a (very) special case of a result by Langton (\cite{Langton}, Proposition 3), who proves, under very general assumptions (for instance if the field extension under consideration is separable and algebraic), that semi-stability is a notion invariant under base change for torsion-free coherent sheaves on a non-singular projective variery.\\ To show that (1) does not necessarily imply (3), we identify all bundles $(\calE,\sigt)$ which are stable in the real (resp. quaternionic) sense. We note that when $\calF$ is any holomorphic vector bundle, there is a commutative diagram

$$\begin{CD}
\os{\cF}  @>\sigt>>  \cF \\
@VVV  @VVV \\
M @>\sigma>> M
\end{CD}$$

\noindent where $\sigt$ is an invertible, $\C$-antilinear map covering $\sig$ and such that $$\sigt \circ \sigt^{-1}=\Id_{\cF},\ \mathrm{and}\ \sigt^{-1} \circ \sigt = \Id_{\os{\cF}}.$$ Therefore, on $\cF \oplus\os{\cF}$, we may define $$\sigt^+ = \begin{pmatrix} 0 & \sigt \\ \sigt^{-1} & 0 \end{pmatrix}\quad \mathrm{and}\quad \sigt^- = \begin{pmatrix} 0 & -\sigt \\ \sigt^{-1} & 0 \end{pmatrix}\, .$$ $\sigt^+$ and $\sigt^-$ are $\C$-antilinear maps from $\cF\oplus\os{\cF}$ to itself, covering $\sig$, and satisfying $$\sigt^+ \circ \sigt^+ = \begin{pmatrix} \Id_{\cF} & 0 \\ 0 & \Id_{\os{\cF}} \end{pmatrix} = \Id_{\cF\oplus\os{\cF}} ,$$ and $$\sigt^- \circ \sigt^- = \begin{pmatrix} -\Id_{\cF} & 0 \\ 0 & -\Id_{\os{\cF}} \end{pmatrix} = -\Id_{\cF\oplus\os{\cF}}.$$ In other words, $(\cF\oplus\os{\cF},\sigt^+)$ is a real bundle, and $(\cF\oplus\os{\cF},\sigt^-)$ is a quaternionic bundle. We also note that, if $(\cE,\sigt)$ is any real (resp. quaternionic) bundle, the bundle $\End(\cE) \simeq \cE^*\otimes\cE$ of endomorphisms  of $\cE$ always has a \textit{real} structure given by $$\xi\otimes v \longmapsto (\ov{\xi\circ \sigt^{-1}}) \otimes \sigt(v).$$ If we still denote $\sigt$ this real structure, the bundle of real (resp. quaternionic) endomorphisms of $(\cE,\sigt)$ is the bundle $\big(\End(\cE)\big)^{\sigt}$ of $\sigt$-invariant elements of $\End(\cE)$.

\begin{proposition}\label{stability_in_the_real_sense}
Let $(\calE,\sigt_{\calE})$ be a stable real (resp. quaternionic) vector bundle.
\begin{enumerate}
\item Then either $(\calE,\sigt_{\calE})$ is geometrically stable, or there exists a holomorphic vector bundle $\calF$, stable in the holomorphic sense, such that $\calE= \calF \oplus \os{\calF}$. In the latter case, if $(\calE,\sigt)$ is real then $\os{\calF}\neq \calF$ and $\sigt_{\calE}=\sigt^+$, and if $(\calE,\sigt)$ is quaternionic, then $\sigt_{\calE} = \sigt^-$.
\item In the geometrically stable case, the set of real (resp. quaternionic) endomorphisms of $(\calE,\sigt_{\cE})$ is $$\big(\End(\cE)\big)^{\sigt_{\cE}} = \{\lambda \Id_{\calE} : \lambda\in\R\} \simeq_{\R} \R,$$ and, if $\calE=\calF \oplus \os{\calF}$, then $$\big(\End(\cE)\big)^{\sigt_{\cE}} = \{(\lambda\Id_{\calF},\ov{\lambda}\Id_{\calF}) : \lambda\in\C\} \simeq_{\R} \C$$.
\end{enumerate}
\end{proposition}

\noindent Note that the isomorphisms given in part (2) of the Proposition are isomorphisms of \textit{real} vector spaces. Also, a real (resp. quaternionic) bundle which is stable in the real (resp. quaternionic) sense but not geometrically stable, is necessarily of even rank.

\begin{proof} Let $(\calE,\sigt_{\calE})$ be a stable real (resp. quaternionic) vector bundle.
\begin{enumerate}
\item Assume that $(\cE,\sigt_{\cE})$ is not geometrically stable. Then there exists a non-trivial sub-bundle $\cF$ of $\cE$ satisfying $\mu(\cF) \geq \mu(\cE)$. Since, by Proposition \ref{semi-stability_and_geometric_semi-stability}, $\cE$ is semi-stable in the holomorphic sense, we in fact have $\mu(\cF)=\mu(\cE)$ and $\cF$ is also semi-stable. As $(\cE,\sigt_{\cE})$ is real (resp. quaternionic), there is a canonical $\C$-linear isomorphism $\phi:\os{\cE} \longrightarrow \cE$ that allows us to identify $\os{\cF}$ with a sub-bundle of $\cE$. We denote $\cE'$ the sub-bundle generated by the $\sigt_{\cE}$-invariant subsheaf $\cF \cap \os{\cF}$ of $\cE$, and $\cE''$ the sub-bundle generated by the $\sigt_{\cE}$-invariant subsheaf $\cF+\os{\cF}$ of $\cE$. Then we have an exact sequence $$0 \longrightarrow \cE' \longrightarrow \cF \oplus \os{\cF} \longrightarrow \cE'' \longrightarrow 0\ ,$$ where the map $\cF\oplus\os{\cF} \longrightarrow \cE''$ is a morphism of real (resp. quaternionic) bundles when $\cF\oplus\os{\cF}$ is endowed with the real structure $\sigt^+$ (resp. the quaternionic structure $\sigt^-$). Assume now that $\cE'$ and $\cE''$ are non-trivial sub-bundles of $\cE$. Since $\cE'$ and $\cE''$ are $\sigt_{\cE}$-invariant sub-bundles of $\cE$ and $\cE$ is stable in the real (resp. quaternionic) sense, one has $$\frac{d'}{r'} := \mu(\cE') < \mu(\cE)\quad \mathrm{and}\quad \frac{d''}{r''} := \mu(\cE'') < \mu(\cE)\ .$$ But $$\mu(\cE) = \mu(\cF) =: \frac{d}{r}$$ so $d'r < dr'$ and $d'' r < d r''$, and therefore $$d'r + d'' r <dr' + d r''.$$ Moreover, since $\deg(\os{\cF})=\deg(\cF)$ and $\rk(\os{\cF}) = \rk(\cF)$, the exact sequence above implies that $d' + d'' = 2d$ and $r' + r'' = 2r$, so $$\frac{d' + d''}{r' + r''} = \frac{2d}{2r} = \frac{d}{r},$$ and therefore $d'r + d'' r = dr' + dr''$, contradicting the strict inequality above. So $\cE' = \{0\}$ and $\cE'' = \cE$, which means that $\cE \simeq \cF \oplus \os{\cF}$ as a real (resp. quaternionic) bundle. The bundle $\cF$ necessarily is stable as a holomorphic bundle, otherwise a non-trivial sub-bundle $\cF'$ of $\cF$ satisfying $\mu(\cF') \geq \mu(\cF)$ gives a non-trivial, $\sigt^{\pm}$-invariant  sub-bundle $\cF'\oplus\os{\cF'}$ of $(\cE,\sigt_{\cE})$ with slope equal to $\mu(\cF') \geq \mu(\cF) =\mu(\cE)$, contradicting the fact that $(\cE,\sigt_{\cE})$ is stable as a real (resp. quaternionic) bundle (note that $\sigt|_{\os{\cF'}}$ maps $\os{\cF'}$ to $\cF'$ by definition of $\os{\cF'}$). Moreover, when $(\cE,\sigt_{\cE})$ is real, $\os{\cF}$ is not isomorphic to $\cF$, otherwise the diagonal embedding $\cF \longrightarrow \cF \oplus \os{\cF} \simeq \cF \oplus\cF$ would provide a $\sigt^+$-invariant sub-bundle, contradicting the stability of $\cE$ as a real bundle. We note that, in the quaternionic case, the diagonal embedding does not provide a $\sigt^-$-invariant sub-bundle and so does not contradict the stability of $\cE$ as a quaternionic bundle. Indeed, we now give an example of a stable quaternionic bundle of the form $(\cF\oplus\os{\cF},\sigt^-)$ with $\cF$ stable as a holomorphic bundle and satisfying $\os{\cF} \simeq \cF$~: consider a real line bundle $(\cL,\sigt)$ on a real algebraic curve $(M,\sig)$ satisfying $M^{\sig}\neq\emptyset$, then $\os{\cL}\simeq \cL$ and $(\cL\oplus\cL,\sigt^-)$ is a stable quaternionic bundle, for a sub-bundle contradicting this would be a quaternionic line bundle on $(M,\sig)$ and there are no quaternionic line bundles on $(M,\sig)$ when $M^{\sig}\neq\emptyset$.
\item If $(\cE,\sigt_{\cE})$ is geometrically stable, then $$\End(\cE) = \{\lambda\Id_{\cE} : \lambda\in \C\} \simeq \C,$$ and the real structure of $\End(\cE)$ acts as $\lambda \longmapsto \ov{\lambda}$ on such endormophisms, so $$\big(\End(\cE)\big)^{\sigt_{\cE}} = \{\lambda\Id_{\cE} : \lambda\in\R\} \simeq_{\R} \R.$$ If $(\cE,\sigt_{\cE})$ is stable but not geometrically stable, then $\cE = \cF \oplus \os{\cF}$ for some $\cF$ stable in the holomorphic sense (so $\os{\cF}$ is also stable in the holomorphic sense), and $$\End(\cE) = \{(\lambda\Id_{\cE}, \mu\Id_{\cE}) : (\lambda,\mu)\in \C \oplus \C\} \simeq \C \oplus \C.$$ The real structure of $\End(\cE)$ acts as $(\lambda,\mu) \longmapsto (\ov{\mu},\ov{\lambda})$ on such endomorphisms, so $$\big(\End(\cF\oplus\os{\cF})\big)^{\sigt_{\cE}}= 
\{(\lambda,\ov{\lambda}) : \lambda\in\C\} \simeq_{\R} \C.$$
\end{enumerate}
\end{proof}

\noindent Proposition \ref{stability_in_the_real_sense} also proves that a bundle $\cE$ which admits a stable but not geometrically stable real structure $\sigt^+$, also admits the stable quaternionic structure $\sigt^-$. The example in the last part of the proof of (1) shows that the converse is not necessarily true. As a final observation, we point out that, when $r\wedge d=1$, a bundle $(\cE,\sigt)$ which is stable in the real or quaternionic sense, necessarily is geometrically stable (as it is geometrically semi-stable, which implies that it is geometrically stable when $r\wedge d=1$). The previous results suggest that, if we want to think of real points of $\Mod$ as moduli of real and quaternionic bundles, we should restrict our attention, either to the case where $r\wedge d=1$, or to the open sub-scheme $\Mods$, whose complex points are isomorphism classes of geometrically stable bundles. The next result formalises this point of view.

\begin{proposition}\label{real_moduli}
Let $(\calE,\sigt)$ and $(\calE',\sigt')$ be two geometrically stable real (resp. quaternionic) bundles, and assume that $\calE$ and $\calE'$ are isomorphic as holomorphic vector bundles. Then $\calE$ and $\calE'$ are isomorphic as real (resp. quaternionic) vector bundles.
\end{proposition}

\begin{proof}
The assumption of the Proposition is that $\phi: \calE'\overset{\simeq}{\longrightarrow} \calE$. Replacing $\sigt'$ with $\varphi\sigt'\varphi^{-1}$ if necessary, we may assume that $\sigt$ and $\sigt'$ are two distinct real structures on the same vector bundle $\calE$. Then $\sigt\sigt'$ is $\C$-linear and, as $\calE$ is stable, this implies that $\sigt\sigt'=\lambda\in\C^*$. This in turn implies that $$\sigt = \lambda (\sigt')^{-1} = \pm\lambda(\sigt') = \pm \lambda \sigt^{-1}(\lambda\cdot) = \pm \lambda (\pm\sigt)(\lambda\cdot) = |\lambda|^2 \sigt,$$ so $\lambda = e^{i\theta}$ for some $\theta \in \R$, whence one obtains $$\sigt = e^{i\theta} \sigt' = e^{i\frac{\theta}{2}} \sigt' (e^{-i\frac{\theta}{2}}\cdot),$$ showing that $\sigt$ and $\sigt'$ are conjugate by an automorphism of $\calE$.
\end{proof}

\noindent Thus, not only do real points of $\Mods$ represent isomorphism classes of geometrically stable bundles that admit either a real or a quaternionic structure (by Proposition \ref{stable_and_Galois_invariant}), but in addition is that real or quaternionic structure unique up to real or quaternionic isomorphism. Moreover, the automorphism group of a geometrically stable real or quaternionic bundle is equal to $\R^*$ by Proposition \ref{stability_in_the_real_sense}. So we see that $\Mods(\R)$ has many of the good properties that one might expect from a coarse moduli space for objects defined over the field of real numbers.

\subsection{Jordan-H\"older filtrations}

Seshadri has shown (\cite{Seshadri}) that, if $\cE$ is a semi-stable holomorphic bundle, it admits a holomorphic Jordan-H\"older filtration $$\{0\} \subset \cE_0 \subset \cE_1 \subset \cdots \subset \cE_l =\cE\, ,$$ the successive quotients of which are stable bundles of slope $\mu(\cE)$. The associated graded object $$\gr(\cE) =  \cE_1/\cE_0 \oplus \cdots \oplus \cE_l/\cE_{l-1}$$ is a direct sum of stable bundles of equal slope and is called a \textbf{poly-stable} bundle (necessarily, the slope of such a direct sum is equal to the slope of any of its terms). Its graded isomorphism class does not depend on the choice of the filtration, and is called the $S$-equivalence class of $\calE$. In this subsection, we analyse the corresponding situation for semi-stable real and quaternionic bundles. We begin with a definition.

\begin{definition}
Let $(\cE,\sigt)$ be a real (resp. quaternionic) bundle. A \textbf{real (resp. quaternionic) Jordan-H\"older filtration}  of $(\cE,\sigt)$ is a filtration $$\{0\} = \cE_0 \subset \cE_1 \subset \cdots \subset \cE_k=\cE$$ by $\sigt$-invariant holomorphic sub-bundles, whose successive quotients are stable in the real (resp. quaternionic) sense.
\end{definition}

\noindent Let us now study Jordan-H\"older filtrations of semi-stable real and quaternionic bundles of fixed slope $\mu$. We denote $\Bun$ the category of semi-stable holomorphic bundles of slope $\mu$. It is an Abelian category. In particular, if $u:\cE_1 \longrightarrow \cE_2$ is a morphism of semi-stable bundles of slope $\mu$, $\Ker u$ and $\Im u$ are semi-stable bundles of slope $\mu$ and there is an isomorphism $\cE/\Ker u \simeq \Im u$. Moreover, $\Bun$ is Artinian, Noetheriean, stable by extensions, and the simple objects of $\Bun$ are the stable bundles of slope $\mu$, which in particular implies the existence of Jordan-H\"older filtrations in the holomorphic sense for semi-stable bundles of slope $\mu$ (\cite{Seshadri, VLP}).

\begin{theorem}\label{Abelian_cat}
Let $\BunR$ (resp. $\BunH$) denote the category of semi-stable real (resp. quaternionic) bundles of slope $\mu$ on $(M,\sig)$. By Proposition \ref{semi-stability_and_geometric_semi-stability}, it is a strict sub-category of the category $\Bun$ of semi-stable holomorphic bundles of slope $\mu$. Moreover~:
\begin{enumerate}
\item If $u:(\cE_1,\sigt_1) \longrightarrow (\cE_2,\sigt_2)$ is a morphism of real (resp. quaternionic) bundles, then the bundles $\Ker u$ and $\Im u$ are semi-stable real (resp. quaternionic) bundles of slope $\mu$, and the isomorphism $\cE/\Ker u \simeq \Im u$ is an isomorphism of real (resp. quaternionic) bundles. As a consequence, $\BunR$ (resp. $\BunH$) is an Abelian category.
\item The Abelian category $\BunR$ (resp. $\BunH$) is Artinian, Noetherian, and stable by extensions. If $(\cE,\sigt)$ is stable in the real (resp. quaternionic) sense, then its endomorphism ring $(\End\, \cE)^{\sigt}$ is a field which is an algebraic extension of $\R$, so it is either $\R$ or $\C$.
\item The simple objects of $\BunR$ (resp. $\BunH$) are the real (resp. quaternionic) bundles of slope $\mu$ on $(M,\sig)$ that are stable in the real (resp. quaternionic) sense. In particular, a semi-stable real (resp. quaternionic) bundle $(\cE,\sigt)$ admits a real (resp. quaternionic) Jordan-H\"older filtration.
\end{enumerate}
\end{theorem}

\begin{proof}~
\begin{enumerate}
\item Since $\BunR$ (resp. $\BunH$) is a sub-category of the Abelian category $\Bun$, it suffices to prove that, if $u:(\cE_1,\sigt_1) \longrightarrow (\cE_2,\sigt_2)$ is a morphism of real (resp. quaternionic) bundles, then the semi-stable bundles of slope $\mu$, $\Ker u$ and $\Im u$, are in fact real (resp. quaternionic) bundles, so they are objects of $\BunR$ (resp. $\BunH$). This follows from the fact that $\Ker u$ is $\sigt_1$-invariant and $\Im u$ is $\sigt_2$-invariant.
\item Because the rank of a vector bundle is finite, it is obvious that decreasing and increasing sequences of sub-bundles are stationary. Moreover, it follows from (1) that, if $u:(\cE_1,\sigt_1) \longrightarrow (\cE_2,\sigt_2)$ is a non-zero morphism between stable real (resp. quaternionic) bundles of equal slope, then $u$ is an isomorphism. In particular, $(\End\, \cE)^{\sigt}$ is a field, which contains $\R$ as the sub-field of scalar endomorphisms. Since the characteristic polynomial of an element in $(\End\, \calE)^{\sigt}$ has real coefficients, the Cayley-Hamilton Theorem implies that the elements of the field $(\End\, \calE)^{\sigt}$ are algebraic over $\R$.
\item Let $(\cE,\sigt)$ be a stable real (resp. quaternionic) bundle of slope $\mu$. Then it does not admit a non-trivial sub-object in $\BunR$ (resp. $\BunH$), for such a sub-object would have slope $\mu$, contradicting the fact that $\cE$ is stable in the real (resp. quaternionic) sense. So $(\cE,\sigt)$ is a simple object in $\BunR$ (resp. $\BunH$). Conversely, if $(\cE,\sigt)$ is a simple object in $\BunR$ (resp. $\BunH$) and $\cF$ is a non-trivial $\sigt$-invariant sub-bundle of $\cE$, then $\mu(\cF) < \mu(\cE)$, because $\mu(\cF \leq \mu(\cE)$ by the semi-stability of $\cE$ and $\mu(\cF) \neq \mu(\cE)$ by the simplicity of $\cE$. So $(\cE,\sigt)$ is in fact stable in the real (resp. quaternionic) sense. The existence of a real (resp. quaternionic) Jordan-H\"older filtration is then proved in the usual way~: since increasing sequences are stationary, there is a strict sub-object $\calF$ of $(\cE,\sigt)$ which is not contained in any strict sub-object. This $\calF$ in turn contains such a strict sub-object, and one constructs in this way a decreasing sequence of sub-objects of $\calE$. As this sequence is stationary, we get a filtration, whose successives quotients are simple by construction (of course, in this particular category, there is a somewhat simpler proof by induction on the rank).
\end{enumerate}
\end{proof}

\noindent The point to make here is that we need to include the real (resp. quaternionic) bundles which are stable but not necessarily geometrically stable in order to guarantee the existence of real (resp. quaternionic) Jordan-H\"older filtrations for semi-stable real (resp. quaternionic) bundles~: a simple real (resp. quaternionic) bundle might only be stable in the real (resp. quaternionic) sense and not geometrically stable, so a semi-stable real (resp. quaternionic) bundle $(\cE,\sigt_{\cE})$ might only admit a Jordan-H\"older filtration whose successive quotients are stable in the real (resp. quaternionic) sense. As an example, consider the bundle $(\cE,\sigt_{\cE}) \simeq (\cF\oplus\os{\cF},\sig^+)$, with $\cF$ stable as holomorphic bundle and such that $\os{\cF} \not\simeq \cF$. Then $(\cE,\sigt_{\cE})$ is stable as a real bundle, so it admits a real Jordan-H\"older filtration of length one, while it admits no Jordan-H\"older filtration the successive quotients of which are geometrically stable real bundles (note that the diagonal embedding is not $\C$-linear when $\os{\cF}\not\simeq\cF$). Moreover, any holomorphic Jordan-H\"older filtration of $\cE\simeq \cF\oplus\os{\cF}$ has length two, showing that it does not coincide with the real Jordan-H\"older filtration in general. 

\noindent The graded object associated to a real (resp. quaternionic) Jordan-H\"older filtration of a semi-stable real (resp. quaternionic) bundle $(\cE,\sigt)$ is a poly-stable object in the sense of the following definition.

\begin{definition}[Poly-stable real and quaternionic bundles]\label{poly-stability}
A real (resp. quaternionic) vector bundle $(\cE,\sigt)$ on $(M,\sig)$ is called \textbf{poly-stable} if there exist real (resp. quaternionic) bundles $(\cF_j,\sigt_j)_{j}$ of equal slope, stable in the real (resp. quaternionic) sense, such that $$\cE \simeq \cF_1 \oplus\cdots \oplus \cF_k$$ and $$\sigt = \sigt_1 \oplus \cdots \oplus \sigt_k.$$ 
\end{definition}

\noindent By Proposition \ref{stability_in_the_real_sense}, a poly-stable real (resp. quaternionic) bundle is poly-stable in the holomorphic sense. We recall that the holomorphic $S$-equivalence class of a semi-stable holomorphic bundle $\cE$ is, by definition (\cite{Seshadri}), the graded isomorphism class of the poly-stable bundle $\gr(\cE)$ associated to any Jordan-H\"older filtration of $\cE$.

\begin{corollary}\label{empty_intersection}
The $S$-equivalence class, as a holomorphic bundle, of a semi-stable real (resp. quaternionic) bundle $(\cE,\sigt)$ contains a poly-stable real (resp. quaternionic) bundle in the sense of Definition \ref{poly-stability}. Any two such objects are isomorphic as real (resp. quaternionic) poly-stable bundles.
\end{corollary}

\noindent In particular, there is a well-defined notion of real (resp. quaternionic) $S$-equivalence class for a semi-stable real (resp. quaternionic) bundle $(\cE,\sigt)$.

\begin{definition}[Real and quaternionic $S$-equivalence classes]
The graded isomorphism class, in the real (resp. quaternionic) sense, of the poly-stable real (resp. quaternionic) bundle $\gr(\cE,\sigt)$ associated to any real (resp. quaternionic) Jordan-H\"older filtration of $(\cE,\sigt)$, is called the \textbf{real (resp. quaternionic) $S$-equivalence class} of $(\cE,\sigt)$.
\end{definition}
 
\begin{proof}[Proof of Corollary \ref{empty_intersection}]
The first part follows from the existence of a real (resp. quaternionic) Jordan-H\"older filtration in the sense of Theorem \ref{Abelian_cat}. As for the second part, it is enough to show that two real (resp. quaternionic) bundles $(\cE_1,\sigt_1)$ and $(\cE_2,\sigt_2)$ which are stable in the real (resp. quaternionic) sense and isomorphic as holomorphic bundles, are in fact isomorphic as real (resp. quaternionic) bundles. Because the holomorphic Jordan-H\"older filtrations of $\cE_1$ and $\cE_2$ must have equal lengths, there are exactly two cases to consider before proceeding by induction~: \
\begin{itemize}
\item $(\cE_1,\sigt_1) \simeq (\cF_1\oplus \os{\cF_1}, \sigt^{\pm})$ and $(\cE_2,\sigt_2) \simeq (\cF_2\oplus \os{\cF_2}, \sigt^{\pm})$, with $\cF_i$ geometrically stable (and not isomorphic to $\os{\cF_i}$ in the real case),
\item $\cE_1$ and $\cE_2$ are geometrically stable.
\end{itemize}
In the first case, the existence of an isomorphism of real (resp. quaternionic) bundles between $(\cE_1,\sigt_1)$ and $(\cE_2,\sigt_2)$ is immediate because, since $\cE_1$ and $\cE_2$ are poly-stable and isomorphic as holomorphic bundles, one has $\cF_1\simeq\cF_2$ or $\cF_1\simeq\os{\cF_2}$. In the second case, the assumption is that there is an isomorphism $\phi: \cE_2\overset{\simeq}{\longrightarrow} \cE_1$ of geometrically stable holomorphic bundles, and this is exactly the situation of Proposition \ref{real_moduli}.
\end{proof}

\noindent We point out that a same poly-stable object may admit, however, both a real and a quaternionic structure, showing that it belongs both to a real and to a quaternionic $S$-equivalence class (for instance, $\cF\oplus\os{\cF}$ admits the real structure $\sigt^+$ and the quaternionic structure $\sigt^-$). A final instructive example is given as follows. 

\begin{example}\label{example}
Let $(\cL,\sigt)$ be a real (resp. quaternionic) line bundle on $(M,\sig)$. Then $\cL\oplus\cL$ admits two non-conjugate, non-stable, real (resp. quaternionic) structures, namely $$\sigt \oplus \sigt \quad \mathrm{and}\quad \sigt^+ = \begin{pmatrix} 0 & \sigt \\ \sigt & 0 \end{pmatrix}.$$ We note that $\sigt^+$ is indeed quaternionic when $\sigt$ is quaternionic. The two non-conjugate poly-stable real (resp. quaternionic) structures $\sigt\oplus\sigt$ and $\sigt^+$ are, however, $S$-equivalent in the real (resp. quaternionic) sense. Indeed, $$\gr(\cL\oplus\cL, \sigt \oplus\sigt) = (\cL,\sigt) \oplus (\cL,\sigt)$$ and $(\cL \oplus \cL, \sigt^+)$ admits the real (resp. quaternionic) Jordan-H\"older filtration $$\{0\} \subset \cL_{\Delta} \subset \cL \oplus \cL\, ,$$ where $\cL_{\Delta}$ is the image of the diagonal embedding
 \begin{eqnarray*} 
 \cL & \longrightarrow & \cL \oplus \cL \\
u & \longmapsto & (u,u) \, .
\end{eqnarray*} 
In particular, $(\cL_{\Delta},\sigt^+|_{\cL_{\Delta}})$ is isomorphic to $(\cL,\sigt)$ as a real (resp. quaternionic) bundle. Moreover, the map 
$$ \begin{array}{ccc} 
 (\cL \oplus \cL) / \cL_{\Delta} & \longrightarrow & \cL \\
(v,w) & \longmapsto & i(v - w)
\end{array}$$ 
is an isomorphism of real (resp. quaternionic) bundles with respect to $\sigt^+$ and $\sigt$, so $$\gr(\cL\oplus\cL,\sigt^+) \simeq (\cL,\sigt) \oplus (\cL,\sigt).$$
\end{example}

\noindent We conclude the present subsection by pointing out that the occurrence of stable objects which are direct sums of stable holomorphic bundles (here, $\cF \oplus \os{\cF}$ with $\cF$ geometrically stable) has appeared before in related contexts, such as the study of orthogonal and spin bundles on a curve (\cite{Ramanan}), or, more recently, the study of $\mathbf{U}(p,q)$-Higgs bundles on a curve (\cite{BGPG}).

\subsection{Real points of moduli schemes of semi-stable vector bundles}\label{subsection:real_points}

We proposed, in subsection \ref{stability}, a notion of moduli (real points of $\Mods$) for geometrically stable real and quaternionic bundles, which includes the usual "good case" where $r\wedge d=1$ (the coprime case)~: in this  case, $\Mod(\C)=\Mods(\C)$, and points of $\Mod(\R)=\Mods(\R)$ are in bijection with isomorphism classes of geometrically stable real and quaternionic bundles of rank $r$ and degree $d$. But when $r$ and $d$ are not coprime, there are holomorphic vector bundles on $X(\C)$ which are semi-stable but not stable, and $\Mod(\C)$ is defined as the set of graded isomorphism class of poly-stable bundles of rank $r$ and degree $d$. As the functor $\cE \longmapsto \os{\cE}$ preserves the rank and degree of a holomorphic vector bundle, it sends a holomorphic Jordan-H\"older filtration $\{0\}=\cE_0 \subset \cE_1 \subset \cdots \subset \cE_l=\cE$ of $\cE$ to the holomorphic Jordan-H\"older filtration $$\{0\}=\os{\cE_0} \subset \os{\cE_1} \subset \cdots \subset \os{\cE_l} =\os{\cE}$$ of $\os{\cE}$, so it induces an involution $$[\gr(\cE)] \longmapsto [\gr(\os{\cE})]$$ of the moduli variety $\Mod(\C)$, which is precisely the $\Gal(\C/\R)$-action induced by the real structure of $X$. If $\cE$ is stable as a holomorphic bundle, then $\gr(\cE) \simeq \cE$, so, if $\cE$ is a real point of $\Mods\subset \Mod$, then Proposition \ref{stable_and_Galois_invariant} shows that $\cE$ is either real or quaternionic, and cannot be both. Moreover, by Proposition \ref{real_moduli}, the real (resp. quaternionic) structure thus defined on $\cE$ is unique up to real (resp. quaternionic) isomorphism. The situation is not quite as nice, however, when $\cE$ is semi-stable but not stable, and satisfies $\gr(\os{\cE}) \simeq \gr(\cE)$. Let us for instance analyse the case where the holomorphic Jordan-H\"older filtration of $\cE$ has length two. Then we write $$\gr(\cE) = \cF_1 \oplus \cF_2\, ,$$ where $\cF_1$ and $\cF_2$ are stable holomorphic bundles of equal slope. So $\os{\cF_1}\oplus\os{\cF_2}$ is isomorphic to $\cF_1 \oplus \cF_2$ if and only if one of the following two options occurs~:
\begin{itemize}
\item[-] $\os{\cF_1}\simeq \cF_1$ and $\os{\cF_2} \simeq \cF_2$. This implies that each $\cF_i$ is either real or quaternionic, and cannot be both. Their direct sum, however, might be of neither type, for instance if $(\cF_1,\sigt_1)$ is real, $(\cF_2,\sigt_2)$ is quaternionic, and $\cF_1 \oplus \cF_2$ is endowed with the $\C$-antilinear map $\sigt_1 \oplus \sigt_2$.
\item[-] $\os{\cF_1} \simeq \cF_2$ and $\os{\cF_2} \simeq \cF_1$. Then $\cE=\cF \oplus \os{\cF}$, and this semi-stable holomorphic bundle may be endowed with the real structure $$\sigt_{\cE}^+ = \begin{pmatrix} 0 & \sigt \\ \sigt^{-1} & 0 \end{pmatrix},$$ or the quaternionic structure $$\sigt_{\cE}^- = \begin{pmatrix} 0 & -\sigt \\ \sigt^{-1} & 0 \end{pmatrix},$$ where $\sigt$ is the invertible $\C$-antilinear map

$$\begin{CD}
\os{\cF}  @>\sigt>>  \cF \\
@VVV  @VVV \\
M @>\sigma>> M.
\end{CD}$$

\end{itemize}
\vskip10pt
\noindent So, in sum, $S$-equivalence classes of semi-stable real or quaternionic bundles always are fixed points of the involution $[\gr(\cE)] \longmapsto [\gr(\os{\cE})]$, but the converse is not true~: there may be fixed points which are $S$-equivalence classes of semi-stable holomorphic bundles that are neither real nor quaternionic (for instance the direct sum of a stable real bundle and a stable quaternionic bundle of equal slope), and there may be fixed points which are $S$-equivalence classes of semi-stable holomorphic bundles that admit both real and quaternionic structures. So we see that real points of $\Mod$ only give a good notion of moduli for semi-stable real and quaternionic bundles when $r\wedge d=1$, in which case such bundles are in fact stable. For arbitrary $r$ and $d$, a real point of $\Mod$ is not necessarily the real (resp. quaternionic) $S$-equivalence class of a semi-stable real (resp. quaternionic) bundle. The gauge-theoretic construction that we present in the next sƒection, however, gives a nice notion of moduli space for a larger class of real (resp. quaternionic) bundles than just the geometrically stable ones, namely those that are semi-stable and real (resp. quaternionic) \textit{and} have a fixed topological type~: in Theorem \ref{points_of_L_sigt}, we show that real (resp. quaternionic) $S$-equivalence classes of such bundles are in bijection with the points of certain Lagrangian quotients defined using anti-symplectic involutions of the space of unitary connections on a fixed real (resp. quaternionic) Hermitian bundle. As a consequence of the construction, these Lagrangian quotients embed onto connected subsets of real points of $\Mod$ and, if we restrict our attention to geometrically stable bundles, we obtain in this fashion (Theorem \ref{connected_components}) exactly the connected components of $\Mods(\R)$.

\section{The differential geometric approach}

\subsection{The momentum map picture}

We now recall the general framework in which we shall prove our result. It is commonly known as the Atiyah-Bott-Donald\-son momentum map picture, and the foundational, key result in this approach is Donaldson's formulation of the Narasimhan-Seshadri theorem (\cite{NS}, \cite{Don_NS}). If one chooses a Hermitian metric on a fixed smooth complex vector bundle of rank $r$ and degree $d$ on the compact Riemann surface $M$, holomorphic structures on $E$ correspond bijectively to $\cxgaugegp$-orbits of unitary connections on $E$ (we denote $\gaugegp$ the group of unitary automorphisms of $E$ and $\cxgaugegp$ the group of all complex linear automorphisms of $E$, respectively called the \textbf{gauge group} and the \textbf{complex gauge group}). Explicitly, the holomorphic sections of the holomorphic bundle $(E,d_A)$ defined by a unitary connection $A$ are the elements of $\ker d_A^{(0,1)}$, the kernel of the $(0,1)$ part of the covariant derivative $$d_A : \sectionsofE \longrightarrow \oneforms = \onezeroforms \oplus \zerooneforms.$$ The space $\conn$ of all unitary connections on $E$ is an infinite-dimensional affine space, whose group of translations is $\antiHermoneforms$, and, provided that one considers $L^2_1$ connections with $L^2$ curvature instead of $C^{\infty}$ such objects, $\conn$ is a Banach manifold with a K\"ahler structure~: its complex structure is induced by the Hodge star of $M$, and the symplectic form is given by $$\w_A(a,b) = \int_M -\tr (a\wedge b)$$ for all $a,b\in T_A\conn \simeq \antiHermoneforms$. Likewise, if one considers  $L_2^2$ gauge transformations, the gauge group is a Banach Lie group, acting on $\conn$ by $$u(A) = A + (d_A\, u) u^{-1}.$$ As noted by Atiyah and Bott, this action is Hamiltonian, the momentum map being the curvature map $$F: \begin{array}{ccc} \conn & \longrightarrow & \antiHermtwoforms \\ A & \longmapsto & F_A \end{array}.$$ In what follows, we denote $*$ the Hodge star of $M$. In particular, it sends a section of $\mathfrak{u}(E)$ to an element in $\antiHermtwoforms$.

\begin{theorem}[Donaldson, \cite{Don_NS}]
A holomorphic vector bundle $\calE$ of rank $r$ and degree $d$ on $M$ is stable if, and only if, the corresponding $\cxgaugegp$-orbit $O(\calE)$ of unitary connections on $E$ contains an irreducible, minimal Yang-Mills connection, meaning a unitary connection $A$ such that~:
\begin{enumerate}
\item $\Stab_{\cxgaugegp}(A) = \C^*$,
\item $$F_A = *
\begin{pmatrix}
i 2\pi \frac{d}{r} & & \\
& \ddots & \\
& & i 2\pi \frac{d}{r}
\end{pmatrix}.$$
\end{enumerate}
\noindent Moreover, such a connection is unique up to a unitary automorphism of $E$.
\end{theorem}

\noindent Connections satisfying Condition (2) are absolute minima of the Yang-Mills functional $A \mapsto \int_M \| F_A \|^2$ on $E$, and the last part of the theorem says that if $A$ and $A'$ are two irreducible, \textit{minimal Yang-Mills} connections which are $\cxgaugegp$-conjugate, then they are $\gaugegp$-conjugate. Poly-stable bundles of rank $r$ and degree $d$ are seen to be those which admit a unitary connection having curvature $F_A = *i 2\pi \frac{d}{r} \Id_E$, so the set of $S$-equivalence classes of semi-stable vector bundles is in bijection with gauge orbits of minimal Yang-Mills connections. In other words, $$\Mod(\C) \simeq \fibre \big/ \gaugegp,$$ a K\"ahler quotient obtained from the infinite-dimensional manifold $\conn$.

\subsection{Real structures on spaces of unitary connections}

We now come to the heart of this paper~: the existence of a finite family of real structures (anti-symplectic, involutive isometries) of the space $\conn$, along with compatible involutions of the gauge group $\gaugegp$ and the space $\antiHermtwoforms$ (which may be identified to $(Lie(\gaugegp))^*$), all of them determined by the choice of a real or quaternionic Hermitian structure $\sigt$ on $E$, and such that the associated Lagrangian quotients $$\calL_{\sigt} := \left( \fibre \right)^{\sigt} \big/ \gaugegp^{\sigt}$$ are the connected components of $\Mod(\R)$ when $r\wedge d =1$. The key, albeit easy, property of the various involutions that we shall consider, is that their fixed points are precisely the unitary connections that define real or quaternionic \textit{holomorphic} structures on the fixed real or quaternionic Hermitian bundle $(E,\sigt)$. It should be noted that it is also an involution that shall characterise unitary connections defining quaternionic holomorphic structures, as opposed to an automorphism squaring to minus the identity.\\ Let $(E,\sigt)$ be a real or quaternionic Hermitian bundle of rank $r$ and degree $d$ ($\sigt^2= + \Id_E$ or $-\Id_E$), and let $\phi: \os{E} \to E$ be the corresponding isomorphism ($\os{\phi} = \phi^{-1}$ or $\os{\phi} = -\phi^{-1}$). For any unitary connection $A\in \conn$ , we define $\ov{A}$ to be the connection given by $$d_{\ov{A}} s = \phi \big( d_{\os{A}} (\phi^{-1} s) \big)$$ for any smooth global section $s$ of $E$ (observe that this formula is similar to the gauge action, except that $\phi$ is not an endomorphism of $E$, and $\os{A}$ is the connection induced by $A$ on $\os{E}$). Thus, starting from a map $$d_A: \sectionsofE \longrightarrow \oneforms,$$ we obtain a new map $$d_{\ov{A}}: \sectionsofE \longrightarrow \oneforms$$ which, locally, is given by $$d_{\ov{A}}s = ds + \phi \os{A} \phi^{-1}.$$ This transformation is involutive, for

\begin{eqnarray*}
\phi\big(\os{(\phi \os{A} \phi^{-1})} \big) \phi^{-1} & = & (\phi \os{\phi}) A (\os{\phi}\,^{-1}\, \phi^{-1}) \\
& = & (\pm \Id_E) A (\pm \Id_E) \\
& = & A,
\end{eqnarray*}

\noindent as the minus signs cancel out in the quaternionic case.

\begin{proposition}
The involution $A \mapsto \ov{A}$ is an anti-symplectic isometry of $\conn$.
\end{proposition}

\begin{proof}
The tangent map to $A\mapsto \ov{A}$ is $$\begin{array}{ccc}
\antiHermoneforms & \longrightarrow & \antiHermoneforms \\
a & \longmapsto & \phi \os{a} \phi^{-1},
\end{array}$$

\noindent which is $\C$-antilinear. Moreover, as $\tr(XY)$ is real-valued on anti-Hermitian matrices and $\sig$ is an orientation-reversing isometry of $M$, one has~:

\begin{eqnarray*}
\w_{\ov{A}} (\ov{a}, \ov{b}) & = & \int_M -\tr(\ov{a} \wedge \ov{b}) \\
& = & \int_M \sig^* \big(-\tr (a \wedge b) \big) \\
& = & - \int_M -\tr(a \wedge b) \\
& = & - \w_A(a,b).
\end{eqnarray*}
\end{proof}

\noindent We now observe that, when $E$ has a real structure, so do all the spaces $\kforms$, and when $E$ has a quaternionic structure, so do all the spaces $\kforms$. In both cases, these structures may be written, in slightly abusive notation (pulling back the differential form followed by applying the real or quaternionic structure), 

$$
\begin{array}{ccc}
\kforms & \longrightarrow & \kforms \\
\eta & \longmapsto & \sigt \circ \eta \circ \sig.
\end{array}
$$

\begin{definition}[Real and quaternionic connections]
Let $(E,\sigt)$ be a real Hermitian bundle. A unitary connection $$d_A : \sectionsofE \longrightarrow \oneforms $$ is called \textbf{real} if it commutes to the respective real structures of $\sectionsofE$ and $\oneforms $. Likewise, a unitary connection on a quaternionic Hermitian bundle $(E,\sigt)$ is called \textbf{quaternionic} if it commutes to the respective quaternionic structures of $\sectionsofE$ and $\oneforms $.
\end{definition}

\noindent The point of this definition is that, if $d_A$ is real, then the real structure of $\sectionsofE$ leaves $\ker d_A^{(0,1)}$ invariant, so $(E,\sigt,d_A)$ is a real \textit{holomorphic} bundle (its space of holomorphic sections has a real structure). Likewise, if $d_A$ is quaternionic, then the quaternionic structure of $\sectionsofE$ induces a quaternionic structure on $\ker d_A^{(0,1)}$, so $(E,\sigt,d_A)$ is a quaternionic holomorphic bundle. The observation is then as follows.

\begin{proposition}\label{charac_real_or_quat}
Let $(E,\sigt)$ be a real (resp. quaternionic) smooth Hermitian bundle, and let $$A \longmapsto \ov{A}$$ be the involution of $\conn$ associated to $\sigt$. Then a unitary connection $A$ on $E$ is real (resp. quaternionic) if and only if $\ov{A} = A$.
\end{proposition}

\begin{proof}
Let us assume that $\sigt$ is a real structure. The key observation is that, for any $s \in \sectionsofE$, $$d_A (\sigt \circ s \circ \sig) = \sigt (d_{\ov{A}} s) \sig,$$ which follows from the definition of $\ov{A}$. Since $A$ is real if and only if $$d_A (\sigt \circ s \circ \sig) = \sigt (d_A s) \sig,$$ we obtain that $A$ is real if, and only if, $d_{\ov{A}} = d_A$. The exact same proof works if $\sigt$ is a quaternionic structure.
\end{proof}

\noindent It remains to prove that the involution $A \mapsto \ov{A}$ preserves the fibre $\fibre$ of the momentum map, and takes a gauge orbit to a gauge orbit. This is classically proved by showing the existence of an involution of $\gaugegp$, inducing an involution of $(Lie\,\gaugegp)^* \simeq \antiHermtwoforms$, both of which are compatible with $A \mapsto \ov{A}$. We define 

$$
\begin{array}{ccc}
\gaugegp & \longrightarrow & \gaugegp \\
u & \longmapsto & \phi\, \os{u}\, \phi^{-1}
\end{array}
$$

\noindent and

$$
\begin{array}{ccc}
\antiHermtwoforms & \longrightarrow & \antiHermtwoforms \\
R & \longmapsto & \phi\, \os{R}\, \phi^{-1},
\end{array}
$$

\noindent where $\phi$ is, as earlier, the isomorphism $\phi: \os{E} \overset{\simeq}{\longrightarrow} E$ determined by the real or quaternionic structure of $E$. Note that $\os{u}$ and $\os{R}$ define endomorphisms of $\os{E}$, so the proposed formulae make sense. In the following, we simply denote $u \mapsto \ov{u}$ and $R \mapsto \ov{R}$ the involutions above. It is a simple matter to verify that the second involution is induced by the first one under the identification $(Lie\,\gaugegp)^* \simeq \antiHermtwoforms$.

\begin{proposition}
One has the following compatibility relations
\begin{enumerate}
\item between the involution $A \mapsto \ov{A}$ and the gauge action~: $$\ov{u(A)} = \ov{u}\, (\ov{A}),$$
\item and between the involution $A \mapsto \ov{A}$ and the momentum map of the gauge action~: $$F_{\ov{A}} = \ov{F_A}.$$
\end{enumerate}
\end{proposition}

\begin{proof}~
\begin{enumerate}
\item One has

\begin{eqnarray*}
\ov{u(A)} & = & \ov{A + (d_A u) u^{-1}} \\
& = & \ov{A} + (\ov{d_A u}) \, \ov{u^{-1}} \\
& = & \ov{A} + (d_{\ov{A}} \ov{u}) \ov{u}^{-1} \\
& = & \ov{u}\, (\ov{A}).
\end{eqnarray*}

\item For all $s \in \sectionsofE$, one has

\begin{eqnarray*}
d_{\ov{A}} (d_{\ov{A}} s) & = & d_{\ov{A}} \left[ \sigt \big( d_A (\sigt \circ s \circ \sig) \big)\sig \right] \\
& = & \sigt \left[ d_A \big(d_A (\sigt \circ s \circ \sig) \big) \right]\sig.
\end{eqnarray*}

\noindent As the $2$-form $F_A$ is determined by the operator $d_A \circ d_A$, and $F_{\ov{A}}$ by $d_{\ov{A}} \circ d_{\ov{A}}$, the above relation between these operators translates to $$F_{\ov{A}} = \ov{F_A}.$$
\end{enumerate}
\end{proof}

\noindent As $*i2\pi\frac{d}{r}\Id_E \in \antiHermtwoforms$ is a fixed point of $R\mapsto \ov{R}$, compatibility relation (1) shows that the involution $A \mapsto \ov{A}$ induces an involution of $\fibre$, and we denote $$\left(\fibre\right)^{\sigt}$$ the fixed-point set of that involution. It consists of minimal Yang-Mills connections which are fixed by the involution $A \mapsto \ov{A}$ (and so are either real or quaternionic, depending on the type of $\sigt$). Compatibility relation (2) shows that the group $\realgaugegp$ of gauge transformations which commute to the real or quaternionic structure of $E$ (called the \textbf{real} or \textbf{quaternionic gauge group}) acts on the set $(\fibre)^{\sigt}$, and the next result, perhaps our most important observation, shows that the intersection of the $\gaugegp$-orbit of a poly-stable real or quaternionic connection with $\conn^{\, \sigt}$, is a single $\realgaugegp$-orbit.

\begin{proposition}\label{single_real_orbit}
Let $A,A'$ be two connections which satisfy $\ov{A} = A$ and $\ov{A'} = A'$, and assume that $A$ and $A'$ define real (resp. quaternionic) structures which are poly-stable. Then $A$ and $A'$ lie in the same $\gaugegp$-orbit if, and only if, they lie in the same $\realgaugegp$-orbit.
\end{proposition}

\begin{proof}
By Propositions \ref{stability_in_the_real_sense} and \ref{charac_real_or_quat} combined with Donaldson's theorem, a connection $A$ defines a poly-stable real (resp. quaternionic) structure, if and only if it is a direct sum $A = A_1 \oplus \cdots \oplus A_k$ of unitary connections $A_i$ such that each $A_i$ is either an irreducible minimal Yang-Mills connection satisfying $\ov{A_i} = A_i$, or is of the form $B_i \oplus \os{B_i}$, with $B_i$ an irreducible minimal Yang-Mills connection. So it suffices to prove the proposition in each of these two cases.
The 'if' part of the proposition is obvious. To prove the 'only if' part, let us assume that $A' = u(A)$ for some $u \in \gaugegp$. Then $$u(A) = A' = \ov{A'} = \ov{u(A)} = \ov{u}\, (\ov{A}) = \ov{u}\, (A).$$ Let us first treat the case where $A$ and $A'$ are irreducible connections. Then $u^{-1}\, \ov{u} \in \gaugegp \cap \C^*=S^1$, so $u^{-1}\, \ov{u} = e^{i\theta}$ for some $\theta \in \R$. Put then $v = e^{i\frac{\theta}{2}} u$. Then $$v(A) = u(A) = A'$$ and $$\ov{v} = e^{-i\frac{\theta}{2}} \ov{u} = e^{-i\frac{\theta}{2}} e^{i\theta} u = e^{i\frac{\theta}{2}} u = v,$$ so $v\in \realgaugegp$. Consider now the case where $A=B\oplus\os{B}$ and $A'=B'\oplus \os{B'}$. Then $u\in \gaugegp$ such that $u(A)=A'$ is of the form $w\oplus \os{w}$, where $w$ is a unitary gauge transformation. Moreover, the stabiliser of $A=B\oplus \os{B}$ in $\gaugegp$ is isomorphic to $S^1 \times S^1$, so now $u^{-1}\ov{u}= e^{i\theta} \oplus e^{-i\theta}$ for some $\theta\in\R$. Then $$v:= (e^{i\frac{\theta}{2}} w) \oplus (e^{-i\frac{\theta}{2}} \os{w}) $$ satisfies $$v(B\oplus \os{B}) = u(B\oplus \os{B}) = B '\oplus \os{B'},$$ and $$\ov{v} = (e^{-i\frac{\theta}{2}} \oplus e^{i\frac{\theta}{2}}) (\ov{w\oplus\os{w}}) = (e^{-i\frac{\theta}{2}} \oplus e^{i\frac{\theta}{2}}) (e^{i\theta} \oplus e^{-i\theta}) (w\oplus\os{w}) = v,$$ so again $v\in\realgaugegp$.
\end{proof}

\noindent Proposition \ref{single_real_orbit} implies that the map $$\calL_{\sigt} \longrightarrow \fibre / \gaugegp\, ,$$ taking the $\realgaugegp$-orbit of a real (resp. quaternionic) minimal Yang-Mills connection to its $\gaugegp$-orbit, is injective. Moreover, the involution $A \longmapsto \ov{A}$ induces the Galois action $[\cE]\longmapsto [\os{\cE}]$ on $\Mod(\C)$, so the Lagrangian quotient $\calL_{\sigt}$ indeed embeds into $\Mod(\R)$. The next result gives the gauge-theoretic construction of moduli spaces of real and quaternionic vector bundles alluded to in the introduction and at the end of Section \ref{vb_over_real_curves}.

\begin{theorem}\label{points_of_L_sigt} Let $(E,\sigt)$ be a real (resp. quaternionic) smooth Hermitian bundle of rank $r$ and degree $d$.
Then $\realgaugegp$-orbits in $(\fibre)^{\sigt}$ correspond bijectively to real (resp. quaternionic) $S$-equivalence classes of semi-stable real (resp. quaternionic) holomorphic vector bundles that are smoothly isomorphic to $(E,\sigt)$. 
\end{theorem}

\noindent In other words, the points of the Lagrangian quotient $$\cL_{\sigt} = \big(\fibre\big)^{\sigt} / \realgaugegp$$ are in bijection with real (resp. quaternionic) $S$-equivalence classes of semi-stable real (resp. quaternionic) holomorphic vector bundles that are smoothly isomorphic to $(E,\sigt)$, making this Lagrangian quotient a moduli space for such bundles in a reasonable sense.

\begin{proof}
The statement will follow from Propositions \ref{charac_real_or_quat} and \ref{single_real_orbit} combined with Donaldson's theorem. Let $(\cE,\sigt_{\cE})$ be a semi-stable real (resp. quaternionic) holomorphic bundle which is smoothly isomorphic to $(E,\sigt)$. Then, by Theorem \ref{Abelian_cat}, $(\cE,\sigt)$ admits a real (resp. quaternionic) Jordan-H\"older filtration, for which $\gr(\cE)$ is a real (resp. quaternionic) poly-stable bundle of the form $$(\cE_1,\sigt_1) \oplus \cdots \oplus (\cE_l,\sigt_l)$$ with each $(\cE_i,\sigt_i)$ stable in the real (resp. quaternionic) sense, and $(E,\sigt)$ is smoothly isomorphic to the direct sum $$(E_1,\sigt_1) \oplus \cdots \oplus (E_l,\sigt_l)\, ,$$ where $E_i$ is the underlying smooth bundle of $\cE_i$. We want to show that there is associated to such a poly-stable real (resp. quaternionic) bundle a uniquely defined $\realgaugegp$-orbit in $(\fibre)^{\sigt}$, where $\sigt$ has now been smoothly identified with $\sigt_1\oplus \cdots \oplus\sigt_l $. By Proposition \ref{stability_in_the_real_sense}, a stable real (resp. quaternionic) bundle $(\cE,\sigt_{\cE})$ is either geometrically stable, or of the form $\cE = \cF \oplus \os{\cF}$, with $\cF$ stable as a holomorphic bundle. In the latter case, the real (resp. quaternionic) structure on $\cF$ is given by $$\sigt^{+} = \begin{pmatrix} 0 & \sigt \\ \sigt^{-1} & 0 \end{pmatrix}\quad \left(\mathrm{resp.}\ \sigt^{-} = \begin{pmatrix} 0 & -\sigt \\ \sigt^{-1} & 0 \end{pmatrix}\right)\, ,$$ where $\sigt$ is the invertible $\C$-antilinear map

$$\begin{CD}
\os{\cF}  @>\sigt>>  \cF \\
@VVV  @VVV \\
M @>\sigma>> M.
\end{CD}$$

\noindent In the geometrically stable case, the holomorphic structure on $\cE$ is defined, by Donaldson's theorem, by an irreducible, minimal Yang-Mills connection $A$ which, by Proposition \ref{charac_real_or_quat}, satisfies $\ov{A}=A$. When $\cE = \cF \oplus \os{\cF}$, the stable holomorphic structure on $\cF$ is defined, by Donaldson's theorem, by an irreducible, minimal Yang-Mills connection $B$ on $F$ (the underlying smooth bundle of $\cF$), and the holomorphic structure of $\cE$ is defined by the unitary connection $$A := \begin{pmatrix} B & 0 \\ 0 & \sigt^{-1} B \sigt \end{pmatrix} $$ on $E=F\oplus\os{F}$. This connection has curvature $*i2\pi\frac{d}{r}\Id_{E}$, and it satisfies

$$\ov{A} = 
\begin{pmatrix} 0 & \pm\sigt \\ \sigt^{-1}& 0 \end{pmatrix}
\begin{pmatrix} B & 0 \\ 0 & \sigt^{-1} B \sigt \end{pmatrix}
\begin{pmatrix} 0 & \sigt \\ \pm\sigt^{-1}& 0 \end{pmatrix}
= \begin{pmatrix} B & 0 \\ 0 & \sigt^{-1} B \sigt \end{pmatrix}
=  A\, .
$$

\noindent This shows the existence of a surjective map from $(\fibre)^{\sigt}$ to the space of real (resp. quaternionic) $S$-equivalence classes of semi-stable real (resp. quaternionic) holomorphic structures on $(E,\sigt)$. By Donaldson's theorem, two connections $$A,A' \in \fibre$$ define the same $S$-equivalence class if and only if they lie in the same $\gaugegp$-orbit , so, combining with Proposition \ref{single_real_orbit}, two connections $A,A'$ in $(\fibre)^{\sigt}$ define the same $S$-equivalence class (as a holomorphic bundle) if and only if they lie in the same $\realgaugegp$-orbit.
\end{proof}

\noindent The last part in the proof above means that the $\cxgaugegp$-orbit of a real (resp. quaternionic) minimal Yang-Mills connection contains a unique $\realgaugegp$-orbit. This fact is a real analogue of the Kempf-Ness theorem, and it is not necessarily true for arbitrary real (resp. quaternionic) connections.

\subsection[Connected components of the set of real points]{Connected components of the set of real points of the moduli scheme}

Our results so far enable us to give the following description of the connected components of $\Mods(\R)$.

\begin{theorem}\label{connected_components}
Let $E$ be a Hermitian vector bundle of rank $r$ and degree $d$ and denote $\calI$ the set of gauge conjugacy classes of Hermitian real or quaternionic structures on $E$. Then $$\Mods(\R) = \bigsqcup_{[\sigt]\in I} \calL_{\sigt}^s\ ,$$ where $\calL_{\sigt}$ is the Lagrangian quotient $$\left( \fibre \right)^{\sigt} \big/ \realgaugegp$$ constructed in the previous subsection, and $\calL_{\sigt}^s$ is the intersection of $\calL_{\sigt}\subset\Mod(\R)$ with $\Mods(\R)$. Moreover, the $\calL_{\sigt}^s$ are non-empty and connected. In particular, the set of connected components of $\Mods(\R)$ is in bijection with the set of topological types of real and quaternionic structures on $E$.
\end{theorem}

\begin{proof}
By Theorem \ref{points_of_L_sigt}, points of $\calL_{\sigt}^s$ are real (resp. quaternionic) classes of geometrically stable real (resp. quaternionic) holomorphic structures on $(E,\sigt)$, and, by Proposition \ref{stable_and_Galois_invariant}, $$\Mods(\R) = \bigcup_{[\sigt]\in I} \calL_{\sigt}^s.$$ Assume first that $[\cE]\in \calL_{\sigt}^s\cap \calL_{\sigt'}^s$ for distinct $\sigt,\sigt'$. By Proposition \ref{stable_and_Galois_invariant}, $\sigt$ and $\sigt'$ are both real or both quaternionic structures. Moreover, by Proposition \ref{real_moduli}, $(\cE,\sigt)$ and $(\cE,\sigt')$ are isomorphic as real (resp. quaternionic) bundles. This implies that $(E,\sigt)$ and $(E,\sigt')$ are isomorphic as (smooth) Hermitian real (resp. quaternionic) vector bundles. In particular, $\sigt$ and $\sigt'$ are gauge conjugate on $E$ so $\calL_{\sigt}^s = \calL_{\sigt'}^s$, proving that the union of all $\calL_{\sigt}^s$ for $[\sigt]\in I$ is disjoint.\\ Conversely, assume that $\sigt' = \psi \sigt \psi^{-1}$ for some $\psi \in \gaugegp$. Let us denote $\asigt$ and $\asigtp$ the involutions $A \mapsto \ov{A}$ induced by $\sigt$ and $\sigt'$, respectively. We also denote $\phisigt:\os{E} \to E$ and $\phisigtp: \os{E} \to E$ the isomorphisms between $\os{E}$ and $E$ induced by $\sigt$ and $\sigt'$, and $\tausigt$ and $\tausigtp$ the involutions $u \mapsto \ov{u}$ of $\gaugegp$ induced by $\sigt$ and $\sigt'$, respectively. Then $\asigtp = v \asigt v^{-1}$ and $\tausigtp = v \tausigtp v^{-1}$, where $v = \psi \tausigt(\psi^{-1}) \in \gaugegp$. Indeed, $$\phisigtp = \phi_{\psi \sigt \psi^{-1}} = \psi\, \phi_{\sigt}\, \os{\psi}^{-1}:\ \os{E}\overset{\simeq}{\longrightarrow} E,$$ so

\begin{eqnarray*}
\asigtp (A) & = & \phisigtp\, \os{A}\, \phisigt'^{-1} \\
& = & (\psi\, \phisigt\, \os{\psi}^{-1})\, \os{A}\, (\os{\psi}\, \phisigt^{-1}\, \psi^{-1}) \\
& = & \psi\, (\phisigt\, \os{\psi}^{-1}\, \phisigt)\, \phisigt^{-1}\, \os{A}\, \phisigt^{-1}\, (\phisigt\, \os{\psi}\, \phisigt^{-1})\, \psi^{-1} \\
& = & \big(\psi \tausigt(\psi^{-1})\big)\, \asigt(A)\, \big(\psi \tausigt(\psi^{-1})\big)^{-1}.
\end{eqnarray*}

\noindent Moreover, for any $u\in \gaugegp$,

\begin{eqnarray*}
\tausigtp(u) & = & \phisigtp\, \os{u}\, \phisigtp^{-1} \\
& = & (\psi\, \phisigt\, \os{\psi}^{-1}) \, \os{u} \, (\os{\psi}\, \phisigt^{-1}\, \psi^{-1}) \\
& = & \psi\, (\phisigt\, \os{\psi}^{-1}\, \phisigt^{-1})\, \phisigt\,\os{u} \, \phisigt^{-1} (\phisigt\, \os{\psi}\, \phisigt^{-1})\, \psi^{-1} \\
& = & \big(\psi \tausigt(\psi^{-1})\big)\, \tausigt(u)\,  \big(\psi \tausigt(\psi^{-1})\big)^{-1}.
\end{eqnarray*}

\noindent So $\asigtp = v \asigt v^{-1}$ and $\tausigtp = v \tausigt v^{-1}$, where $v\in \gaugegp$ is of the form $\psi \tausigt(\psi^{-1})$ (a symmetric element). This readily implies that the gauge transformation $A \mapsto vAv^{-1}$ establishes a bijection between $\conn^{\asigt}$ and $\conn^{\asigtp}$, which in turn induces a bijection between $\calL_{\sigt}^s$ and $\calL_{\sigt'}^s$. But a gauge transformation on $\conn$ is in fact the identity on $\Mod(\C) = \fibre \big/ \gaugegp$, so $\calL_{\sigt'}^s = \calL_{\sigt}^s$. Therefore, we have proved that $$\Mods(\R) = \bigsqcup_{[\sigt] \in \calI} \calL_{\sigt}^s.$$ It is nice to observe that this is a direct, simple consequence of the general theory of anti-symplectic involutions on Hamiltonian spaces. There remains to study the non-emptiness and the connectedness of the $\calL_{\sigt}^s$. Non-emptiness of $\calL_{\sigt}^s$ is equivalent to the existence of an irreducible, minimal Yang-Mills real (resp. quaternionic) connection on $(E,\sigt)$. But minimal Yang-Mills connections are absolute minima of the Yang-Mills functional and Daskalopoulos has shown that the gradient flow of that functional converges (\cite{Dask}). Moreover, as $$\int_M \| F_{\ov{A}} \|^2 = \int_M \| \ov{F_A} \|^2 = \int_M \| F_A \|^2,$$ the gradient flow of this functional takes a real or quaternionic connection to a real or quaternionic connection, so the limiting connection is of the same type as the original one when following gradient flows. This is used in \cite{BHH}, where Biswas, Hurtubise and Huisman also compute the Morse indices of non-minimal critical sets of the Yang-Mills functional restricted to real and quaternionic connections, to show that, what in our notation is $\calL_{\sigt}^s$, is non-empty and connected (Theorem 6.7 in \cite{BHH}). Note, however, that their point of view is different from ours, as they do not show that their moduli spaces of real and quaternionic vector bundles embed onto connected subsets of real points inside $\Mods(\C)$, nor do they obtain a presentation of them as Lagrangian quotients. As the $\calL_{\sigt}^s$ are connected and closed, with dimension that of $\Mods(\R)$, these are exactly the connected components of $\Mods(\R)$, and counting those connected components amounts to counting the possible topological types of real and quaternionic vector bundles.
\end{proof}

\begin{corollary}
Let $r\wedge d=1$. Then $$\Mod(\R) = \bigsqcup_{[\sigt] \in \calI} \calL_{\sigt}\, .$$ Moreover, the $\calL_{\sigt}$ are non-empty and connected. In particular, the set of connected components of $\Mod(\R)$ when $r\wedge d=1$ is in bijection with the set of topological types of real and quaternionic structures on $E$.
\end{corollary}

\noindent Note that, by Corollary \ref{empty_intersection} and Theorem \ref{points_of_L_sigt}, the intersection, when $r\wedge d\neq 1$, of $\calL_{\sigt}$ and $\calL_{\sigt'}$ for non-conjugate $\sigt,\sigt'$ on the smooth bundle $E$, can only be non-empty if one of them is real and the other one quaternionic (in Example \ref{example}, the \textit{smooth} extension of $(\cL\oplus\cL,\sigt^+)$ introduced there splits, so the underlying Hermitian bundle is indeed isomorphic to that of $(\cL\oplus\cL, \sigt\oplus\sigt)$; in particular, there is indeed only one topological type of quaternionic structure). Moreover, the intersection of an $\cL_{\sigt}$ with an $\cL_{\sigt'}$, if non-empty, is contained in the strictly semi-stable locus of $\Mod(\C)$. We can now prove our main result.

\begin{proof}[Proof of Theorem \ref{main_result}] By Theorem \ref{top_type}, the topological type of a quaternionic vector bundle is determined by its rank and degree, so there is at most one connected component corresponding to quaternionic bundles. When $k(X)=0$ (that is, when $X$ has no real points), there may be at most one connected component of real bundles, and it happens if and only if $d$ is even. Combined with Theorem \ref{top_type}, this gives part (2) of our main result (Theorem \ref{main_result}). To obtain part (1), the part when $k(X) > 0$, it suffices to count the topological types of real bundles, that is, the number of solutions to the equation $$w_1 + w_2 + \cdots + w_{k(X)} = d\ \mod{2}$$ in $\Z / 2\Z$, which is $2^{k(X) -1}$. Precisely, the connected components of $\Mods(\R)$ can be described as follows.

\begin{enumerate}

\item Assume that $k > 0$.

\begin{enumerate}

\item If $r \equiv 1\ (\mod{2})$, then $\Mods(\R)$ has $2^{k-1}$ real components and no quaternionic component (recall that if $k>0$, the rank of a quaternionic bundle must be even). We prove below that any two of these components are homeomorphic.

\item If $r \equiv 0\ (\mod{2})$ and $d \equiv 1\ (\mod{2})$, then $\Mods(\R)$ has $2^{k-1}$ real components and no quaternionic component.

\item If $r \equiv 0\ (\mod{2})$ and $d \equiv 0\ (\mod{2})$, then $\Mods(\R)$ has $2^{k-1}$ real components and $1$ quaternionic component.

\end{enumerate}

\item Assume that $k=0$.

\begin{enumerate}

\item If $r(g-1) \equiv 0\ (\mod{2})$ and $d \equiv 0\ (\mod{2})$, then $\Mods(\R)$ has one real and one quaternionic component. We prove below that these two components are homeomorphic when $g \equiv 1 (\mod{2})$.

\item If $r(g-1) \equiv 0\ (\mod{2})$ and $d \equiv 1\ (\mod{2})$, then $\Mods(\R)$ is empty.

\item If $r(g-1) \equiv 1\ (\mod{2})$, then $\Mods(\R)$ has one real component if $d$ is even, and one quaternionic component if $d$ is odd.

\end{enumerate}

\end{enumerate}

\noindent To finish the proof, it only remains to show that certain identified connected components of $\Mods(\R)$ are homeomorphic.

\begin{enumerate}

\item Assume that $k>0$, and consider a real holomorphic line bundle $(\calL_{\R},\sigt_{\calL_{\R}})$ of degree $0$ on $X(\C)$, with real invariants $(w^{(1)}, \cdots, w^{(k)})\in (\Z/2\Z)^k$ such that $w^{(1)}+\cdots+w^{(k)}=0$. If $(\calE,\sigt)$ is a real bundle on $X(\C)$, of rank $r$ and degree $d$, so is $\calE\otimes \calL_{\R}$, and, over the $j$-th connected component of $X(\R)$, one has $$w_1(\calE_j^{\sigt} \otimes (\calL_{\R})^{\sigt_{\calL_{\R}}}_j) = w_1(\calE_j^{\sigt}) + (r\ \mod{2}) w^{(j)}.$$ In particular, when $r$ is odd, the functor $\calE \mapsto \calE \otimes \calL_{\R}$ induces a homeomorphism between the connected component of $\Mods(\R)$ containing $\calE$ and that containing $\calE \otimes \calL_{\R}$~: by choosing an $\calL_{\R}$ with the appropriate topological invariants $(w^{(1)},\cdots,w^{(k)})$, one may pass from a given connected component of $\Mods(\R)$ to any given other. When $r$ is even, the map induced by the functor $\calE \mapsto \calE \otimes \calL_{\R}$ preserves any given connected component of $\Mod(\R)$, so we do not know if they are homeomorphic. We also note that, when $k>0$, there are no quaternionic line bundles, which prevents one from applying the same technique. 

\item Assume that $k=0$, $r(g-1)$ even and $d$ even. Then there exists a quaternionic line bundle $(\calL_{\H}, \sigt_{\calL_{\H}})$ of degree $0$ on $X(\C)$ if, and only if, $g$ is odd. The functor $\calE \mapsto \calE \otimes \calL_{\H}$ then provides the desired homeomorphism between the real and quaternionic components in that case (note that now $\calE \otimes \calL_{\H}$ is quaternionic, as $\sigt \otimes \sig_{\calL_{\H}}$ squares to $\sigt^2 \otimes \sig_{\calL_{\H}}^2 = \Id_{\calE} \otimes (-\Id_{\calL_{\H}})$). 

\end{enumerate}

\noindent This completes the proof of Theorem \ref{main_result}, thus generalising to rank $r >1$ the results of Gross and Harris on the Picard scheme of a real algebraic curve (\cite{GH}).

\end{proof}

%\bibliographystyle{alpha}
%\bibliography{biblio}

\end{document}